\documentclass[12pt]{article} \textwidth=6in \oddsidemargin=0in
\begin{document}

\def\ch{\raisebox{.4ex}{$\chi$}} \def\st{\tilde{\s}}
\def\chp{\ch^+} \def\chm{\ch^-}
\def\d{\det\,} \def\noi{\noindent} \def\s{\sigma}
\def\t{\tau} \def\W{W_\a} \def\x{\xi} \def\tl{\tilde} \def\tn{\otimes}
\def\ph{\phi} \def\tr{{\rm tr\,}}
\def\iy{\infty} \def\be{\begin{equation}} \def\ee{\end{equation}}
\def\inv{^{-1}} \def\({\left(} \def\){\right)}
\def\ov{\over} \def\a{\alpha}  \def\la{\lambda}
 \def\d{\det\,}  \def\dl{\delta} \def\bR{{\bf R}}\def\D{\Delta} \def\eq{\equiv} \def\e{\eta} \def\z{\zeta} \def\r{\rho} \def\st{\tl \s} \def\rt{\tl \r} \def\rmt{\wt{\r_-}} \def\sp{\vspace{1ex}} \def\noi{\noindent}
 \def\wt{\widetilde} \def\Oda{O(e^{-\dl\a})}  \def\cp{c_+} \def\cm{c_-} \def\up{u_+}\def\um{u_-}  
\def\spt{\wt{\s_+}} \def\smt{\wt{\s_-}} \def\umt{\wt{u_-}}
\def\pmp{{\pm p}} \def\F{\mathcal{F}}

\hfill October 2, 2006\begin{center}{ \large\bf  On the Inverse and Determinant\\ \vskip1ex  of Certain Truncated Wiener-Hopf Operators}\end{center}

\begin{center}{\bf Harold Widom}\end{center}

\begin{center}{Department of Mathematics\\
University of California\\
Santa Cruz, CA 95064}
\end{center}

\begin{center}{\bf 1. Introduction and statement of results}\end{center}

The solution of a problem arising in integrable systems \cite{W1} requires the asymptotics as $\a\to\iy$, with very small error, of the inverses and determinants of truncated Wiener-Hopf operators $\W(\s)$ acting on $L^2(0,\,\a)$, both in the regular case (where the Wiener-Hopf operator $W(\s)$ on $L^2(\bR^+)$ is invertible) and in singular cases. 
This paper treats two cases where $\s$ has simple Fisher-Hartwig singularities, one double zero or two simple zeros. 

We first state a result for the regular case. Assume $\s$ belongs to the Wiener algebra  of functions of the form $a+\hat{k}$ with $k\in L^1(\bR)$,\footnote{The Fourier transform we use is
$\hat{k}(\x)=\int_{-\iy}^\iy e^{ix\x}\,k(x)\,dx$. For the Toeplitz case with a weaker assumption see \cite[Th. 2.14]{BS2}. For (\ref{reginv}) itself see footnote \ref{reg}.} and that
\be\s(\x)\ne0,\ \ \ \Delta\,{\rm arg}\;\s(\x)\Big|_{-\iy}^\iy=0.\label{index}\ee
Then
\be\W(\s)\inv=P_\a\,W(\s)\inv\,P_\a
-Q_\a\,H(1/\smt)\,H(1/\s_+)\,Q_\a+o(1).\label{reginv}\ee 

The notation here is the following. The operators $P_\a$ is multiplication by $\ch_{(0,\a)}$ or extension by zero from $(0,\a)$ to $(0,\iy)$, depending on the context. For functions defined on $(0,\a)$ we define $(Q_\a\,f)(x)=f(\a-x)$, we use $H$ to denote Hankel operators as usual, and for a function $v$ we define $\tl{v}(\x)=v(-\x)$. The term $o(1)$ denotes a family of operators whose operator norms tend to zero as $\a\to\iy$. The functions $\s_\pm$ are the Wiener-Hopf factors of $\s$: their product equals $\s$, and $\s_-^{\pm1}$ resp. $\s_+^{\pm1}$ extend to bounded analytic functions in the lower resp. upper half plane. 

If in addition $\log\s\in L^1$ and  $\int_{-\iy}^\iy|x|\,|k(x)|^2\,dx<\iy$ then the Kac-Achieser formula holds:
\be\d\,\W(\s)\sim G(\s)^\a\,E(\s),\label{regdet}\ee
where
\[G(\s)=\exp\left\{{1\ov2\pi}\int_{-\iy}^{\iy} \log\,\s(\x)\,d\x\right\},\ \ \ 
E(\s)=\exp\left\{{1\ov2\pi i}\int_{-\iy}^{\iy}(\log\s_+(\x))'\,\log\s_-(\x)\,d\x\right\}.\]

In the present paper we assume that 
\[\s(\x)=(\x^2-p^2)\,\t(\x),\]
where $p$ is real, and $\t$ is nonzero, has index zero, and its Wiener-Hopf factors $\t_\pm$ are of the order $\x\inv$ as $\x\to\iy$.

In order to minimize technical details we make a very strong assumption, namely that $e^{\dl\,|x|}k(x)\in L^1(\bR)$ for some $\dl>0$.\footnote{We shall consistently use $\dl$ to denote some positive quantity, different for each occurrence.} This will result in an exponentially small error in the approximation.

We make a choice of which Wiener-Hopf factor of $\s$ incorporates which zero, and take the factors to be
\[\s_-(\x)=i(\x-p)\,\t_-(\x),\quad \s_+(\x)=-i(\x+p)\,\t_+(\x).\]
As functions, we have
\[{\t_-(p)\ov\s_-(\x)}={1\ov i(\x-p)}+u_-(\x),\ \ \ { \t_+(p)\ov\s_+(\x)}={1\ov -i(\x+p)}+u_+(\x),\]
where $u_\pm$ are bounded smooth functions. In fact we think of $\s_\pm(\x)\inv$ as distributions defined by
\[{\t_-(p)\ov\s_-(\x)}={1\ov 0+i(\x-p)}+u_-(\x),\ \ \ { \t_+(p)\ov\s_+(\x)}={1\ov 0-i(\x+p)}+u_+(\x).\]
Observe that if we set $e_p(x)=e^{ipx}$ then $(0-i(\x+p))\inv$ is the Fourier transform of $\chp e_p$ and $(0+i(\x-p))\inv$ is the Fourier transform of $\chm e_p$. (We write $\ch^\pm$ for $\ch_{\bR^\pm}$.)

To state our first result for $p\ne0$ we define
\[\z=H(1/\s_+)\,e_p,\ \ \ \e=H(1/\smt)\,e_p,\]
(here $e_p$ is thought of as restricted to $\bR^+$) and
\[A={2ip\ov\t_-(p)\,\t_+(-p)\,e^{-i\a p}-\t_+(p)\,\t_-(-p)\,e^{i\a p}}.\]
\sp

\noi{\bf Theorem 1} {\boldmath{$(p\ne0)$}}. There is a $\dl>0$ such that if $A=O(e^{\dl\a})$ then $\W(\s)$ is invertible for large $\a$  and
\[\W(\s)\inv\eq P_\a\,W(1/\s_+)\,W(1/\s_-)\,P_\a-Q_\a\,H(1/\smt)\,H(1/\s_+)\,Q_\a\]
\[+A\,[e^{i\a p}\,\t_+(p)\,\t_-(-p)\,P_\a\,\z\tn 
P_\a\,\e-\t_+(p)\,\t_+(-p)\,P_\a\,\z\tn Q_\a\z\]
\[-\t_-(p)\,\t_-(-p)\,Q_\a \e\tn P_\a\,\e+e^{i\a p}\,\t_+(p)\,\t_-(-p)\,Q_\a\e\tn Q_\a\z].\]
\sp

This requires lots of explanation. First, the sign $\eq$ between two operators indicates that the difference is an integral operator whose kernel is uniformly $O(e^{-\dl\a})$ for some $\dl>0$.

There are tensor products in the statement. For functions $v,\,w$ we define $v\tn w$ to be the operator taking a function $f$ to $v\,(w,\,f)$. (The inner product is the integral of the product---no complex conjugate.) In general we have $v_1\tn w_1\;\cdot v_2\tn w_2=(w_1,\,v_2)\,v_1\tn w_2$ and $(v\tn w)\cdot T=v\tn T'w$, where $T'$ is the transpose of $T$.

Next, we define $W((0-i(\x+p))\inv)$ to be convolution by $\chp e_p$ on $(0,\,\iy)$. This is defined on any locally integrable function since
\[((\chp e_p)*f)(x)=e^{ipx}\,\int_0^x e^{-ipy}\,f(y)\,dy.\]
Our assumption on $\s$ implies that $u_+$ is the Fourier transform of an exponentially decaying function,\footnote{Our assumtion implies that $\s_+(\x)$ extends analytically to a strip around the real line and that $\s_+(\x-i\dl)\inv$ is in the Wiener algebra for sufficientlly small $\dl>0$. Therefore, since $(0-i(\x-i\dl))\inv=-(\dl+i\x)\inv$ is in the Wiener algebra, so is $u_+(\x-i\dl)$. This implies that $u_+$ is the Fourier transform of a function with bound $O(e^{-\dl x})$.}
so $W(u_+)$ is defined on any function of at most polynomial growth. This extends the domain of $W(1/\s_+)$ to any function of at most polynomial growth. Similarly we define $W((0+i(\x-p))\inv)$ to be convolution by $\chm e_p$. This is defined only for functions in $L^1$ since
\[((\chm e_p)*f)(x)=e^{ipx}\,\int_x^\iy e^{-ipy}\,f(y)\,dy.\]
Thus $W(1/\s_-)$ acts on $L^1(\bR^+)$. 

Finally, there appear Hankel operators whose symbols are distributions. The symbol of $H(1/\s_+)$ is the Fourier transform of a constant times $\chp e_p$ plus a rapidly decreasing function. Thus the operator has a singular part which can be written as a constant times $e_p\tn e_p$. Similarly for $H(1/\smt)$. In the statement of the theorem we see the product of two of these, or one of them acting on $e_p$. We interpret these by going to the Fourier transforms. For example, from the general identity
\[(v,\,w)={1\ov 2\pi}\,\int_{-\iy}^\iy \hat{v}(-\x)\,\hat{w}(\x)\,d\x\]
we obtain
\[(\chp e_p,\,\chp e_p)= {1\ov 2\pi}\,\int{1\ov 0+i(\x-p)}\,{1\ov 0-i(\x+p)}\,d\x=-{1\ov 2ip},\]
since by going into the upper half-plane we pass the pole of the first factor at $p+0i$. Similarly we find that $H((0-i(\x+p)\inv)\,e_p=-e_p/2ip$. (Everything can be comuputed by first giving $p$ a small positive imaginary part and then passing to the limit.)

Although we shall be using distributions throughout they will all be of a simple form: Fourier transforms of linear combinations of $e_\pmp$ (combinations of 1 and $x$ when $p=0$) and exponentially decaying functions, and their restrictions to $\bR^\pm$. Thus the meaning of the operators $v\to v_\pm$ in Fourier transform space is clear: take the inverse Fourier transform, multiply by $\ch^\pm$, and then take the Fourier transform.

The statement for $p=0$ is simpler since Hankel operators with distribution symbol do not appear. We define
\[B= \a+i\,\({\s_+\ov\s_-}\)'(0),\]
where the factors are chosen so that $\t_-(0)=\t_+(0)$.\footnote{Without this normalization the second term is to be multiplied by $(\t_-/\t_+)(0)=-(\s_-/\s_+)(0)$.\label{B}}
\sp

\noi{\bf Theorem 2} {\boldmath{$(p=0)$}}. The operator $\W(\s)$ is invertible for sufficiently large $\a$ and
\[\W(\s)\inv\eq P_\a\,W(1/\s_+)\,W(1/\s_-)\,P_\a
-Q_\a\,H(\umt)\,H(u_+)\,Q_\a\]
\be-B\inv\,[Q_\a\,H(\umt)1+P_\a\,W(1/\s_+)1]\tn[Q_\a\,H(u_+)1+P_\a\,W(1/\smt)1].\label{singinv}\ee

\sp

There are several results in the literature giving first-order asymptotic results for the inverse when the symbol has a single singularity or zero. The method of \cite{W} gives a first-order result in terms of convergence in operator norm. Rambour and Seghier \cite{RS} considered the Toeplitz case for a symbol with a zero of arbitrary even order and obtained a stronger result in these cases. Ehrhardt \cite{E} obtained results in the Toeplitz case for a general class of symbols with several Fisher-Hartwig singularities, with convergence in the operator norms of weighted $\ell^2$ spaces, although they do not apply to symbols with two real zeros. In this work, too, distributions played an important role. None of these serves our present purpose, which is to find formulas that hold uniformly throughout $(0,\,\a)$ with very small error.

An earlier and more involved derivation of (\ref{singinv}) than the one we give here is in \cite{W2}. The present proof will use a uniformity in the lemma to Theorem~1 and then a passage to the limit. The computation would be tedious if done by hand and we resorted to computer computation. For the skeptical reader we give a direct proof in \S 6, along the lines of the proof of Theorem~1.

For the determinants we assume that in addition
$\log\s\in L^1$. We use $\eq$ here to indicate that the difference of the two sides is $\Oda$ for some $\dl>0$. Now we define
\[G(\s)=\exp\left\{{1\ov2\pi}\int_{-\iy}^{\iy} (\log\s_++\log\s_-)\,d\x\right\},\]
\[E(\s)=\exp\left\{{1\ov2\pi i}\int_{-\iy}^{\iy}(\log\s_+)'\,\log\s_-\,d\x\right\}.\]
In the integrals $\log\s_+$ is interpreted as $\log\,(\x+p+0i)$ plus a regular function while $\log\s_-$ is $\log\,(\x-p-0i)$ plus a regular function. These depend on our choice of factors $\s_\pm$.

For $p\ne0$ the result is
\sp

\noi{\bf Theorem 3} {\boldmath{$(p\ne0)$}}. 
We have as $\a\to\iy$  
\be{\d\W(\s)\ov G(\s)^\a}\eq 
{\t_-(p)\,\t_+(-p)-\t_+(p)\,\t_-(-p)\,e^{2i\a p}\ov\t_+(-p)\,\t_-(p)}\,E(\s).\label{detpne0}\ee
\sp

For $p=0$ the result is
\sp

\noi{\bf Theorem 4} {\boldmath{$(p=0)$}}. We have as $\a\to\iy$
\be{\d\,\W(\s)\ov G(\s)^\a}\eq B\;E(\s).\label{detp=0}\ee
\sp

This  may be rewritten
\[{\d\,\W(\s)\ov G(\s)^\a}=\a\,E(\s)+i\,\({\s_+\ov\s_-}\)'(0)\;E(\s)+\Oda.\]
Thus it gives the second-order asymptotics with very small error. The first-order asymptotics were obtained by Mikaelyan \cite{M}.
In the special case $\s_0(\x)=\x^2/(1+\x^2)$ one has
\[\W(\s_0)=e^{-\a}\,(1+\a/2)\]
exactly, as was also shown in \cite{M}. Exact formulas for arbitrary rational symbols were obtained by B\"ottcher \cite{B} and used by him to obtain further asymptotic results in these cases. For the Toeplitz analogue the asymptotics are known for symbols with any number of zeros \cite{BS1}. (Or see \cite[\S 10.47]{BS}.) For Wiener-Hopf operators the case of two simple zeros was considered in \cite{AM}, but the result was obtained only for $\t$ real and positive. In \S 5 we show that our formula agrees with the one obtained there.

Here is how the approximate inverse in Theorem 1 is obtained. Since (\ref{reginv}) holds for regular symbols we consider 
\be P_\a\,W(\s_+\inv)\,W(\s_-\inv)\,P_\a-Q_\a\,H(1/\smt)\,H(1/\s_+)\,Q_\a\label{approx}\ee
a possible first approximation for singular symbols. So we see what happens when we multiply it by $\W(\s)$ on the left. We find that, with small error, the result is $I$ plus an operator of rank two. The form of this operator suggests the ingredients of a rank two operator that should have been added to (\ref{approx}) to obtain $I$ in the end. A simple computation tells us exactly which operator to take, and the rest is straightforward. For the determinants we introduce a parameter, use the formula for the logarithmic derivative in terms of the inverse, and apply the earlier results.
   
\begin{center}{\bf 2. Proof of Theorem 1}\end{center}

\noi{\bf Lemma}. Define
\[T_\a=P_\a\,W(1/\s_+)\,W(1/\s_-)\,P_\a-Q_\a\,H(1/\smt)\,H(1/\s_+)\,Q_\a\]
\[+A\,\left[e^{i\a p}\,\t_+(p)\,\t_-(-p)\,P_\a\,\z\tn P_\a\,\e-\t_+(p)\,\t_+(-p)\,P_\a\,\z\tn Q_\a\z\right.\]
\be\left.-\t_-(p)\,\t_-(-p)\,Q_\a \e\tn P_\a\,\e+e^{i\a p}\,\t_+(p)\,\t_-(-p)\,Q_\a\z\tn Q_\a\e\right].\label{Ta}\ee
Then if $\dl$ is small enough we have, uniformly in $p$,
\[\W(\s)\,T_\a=I+O(A\,e^{-\dl\a}),\quad T_\a\,\W(\s)=I+O(A\,e^{-\dl\a}).\]

\noi{\bf Proof}. We shall use, innumerable times and without comment, the identities
\[ W(\s_1)\,\W(\s_2)=\W(\s_1\,\s_2)-P_\a\,H(\s_1)\,H(\wt{\s_2})\,P_\a
-Q_\a\,H(\wt{\s_1})\,H(\s_2)\,Q_\a,\]
and
\[ H(\s_1\,\s_2)=W(\s_1)\,H(\s_2)+H(\s_1)\,W(\tl{\s_2}),\] 
and their familiar consequences. (For the Toeplitz analogues of the stated identities see \cite[\S 2.17 and \S7.8]{BS}.)

Let us multiply (\ref{approx}) by $\W(\s)$ on the left. We have
\[\W(\s)\,P_\a\,W(1/\s_+)\,W(1/\s_-)\,P_\a=\W(\s)\,\W(1/\s_+)\,W(\s_-\inv)\,P_\a\]
\[=[\W(\s_-)-Q_a\,H(\st)\,H(1/\s_+)\,Q_\a]\,W(1/\s_-)\,P_\a=I-Q_a\,H(\st)\,H(1/\s_+)\,P_\a\,W(1/\smt)\,Q_\a\]
\[=I-Q_a\,H(\st)\,H(1/\s_+)\,W(1/\smt)\,Q_\a+Q_a\,H(\st)\,H(1/\s_+)\,(I-P_\a)\,W(\smt\inv)\,Q_\a.\]
The first two terms combine as  
\[I-Q_\a\,H(\st)\,H(1/(\s_+\s_-))\,Q_\a.\]
Similarly,
\[-\W(\s)\,Q_\a\,H(1/\smt)\,H(1/\s_+)\,Q_\a=-Q_\a\,W(\st)\,P_\a\,H(1/\smt)\,H(1/\s_+)\,Q_\a\]
\[=-Q_\a\,W(\st)\,H(1/\smt)\,H(1/\s_+)\,Q_\a+
Q_\a\,W(\st)\,(I-P_\a)\,H(1/\smt)\,H(1/\s_+)\,Q_\a.\]
The first term here is
\[Q_\a\,H(\st)\,W(1/\s_-)\,H(1/\s_+)\,Q_\a=Q_\a\,H(\st)\,H(1/(\s_-\,\s_+))\,Q_\a,\]
cancelling the same term above. It follows that $\W(\s)$ times (\ref{approx}) equals $I$ plus
\[Q_a\,H(\st)\,H(1/\s_+)\,(I-P_\a)\,W(1/\smt)\,Q_\a+Q_\a\,W(\st)\,(I-P_\a)\,H(1/\smt)\,H(1/\s_+)\,Q_\a.\footnote{In the regular case with the stated conditions the Hankel operators are all compact, so this operator is $o(1)$. Similarly for the one obtained by multiplying by $\W(\s)$ on the right. This gives (\ref{reginv}) almost immediately. With the stronger assumption on $\s$ this operator is $\Oda$.\label{reg}}\]

For convenience we set
\[c_-=\t_-(p),\quad c_+=\t_+(-p).\]
The Hankel operator $H(1/\s_+)$ on the left above is the sum of two terms. One is $c_+\inv\,e_p\tn e_p$ while the other is $H(u_+)$. Consider the contribution of the latter. It has kernel $O(e^{-\dl\,(x+y)})$, so $H(u_+)\,(I-P_\a)$ has kernel $O(e^{-\dl\,\a}\,e^{-\dl\,(x+y)})$. Since the kernels of both $H(\st)$ and $(I-P_\a)\,W(1/\smt)\,Q_\a$ are bounded it follows that the contribution of the $H(u_+)$ summand of $H(1/\s_+)$ in the first term is $\Oda$. Similarly so is the contribution of the $H(\umt)$ summand of $H(1/\smt)$ in the second term. All this is uniform in $p$. Therefore with this error the above may be replaced by 
\[{1\ov c_+}\,Q_a\,H(\st)e_p\tn Q_\a\,W(1/\s_-)\,(I-P_\a)\,e_p+
{1\ov c_-}\,Q_\a\,W(\st)(I-P_\a)\,e_p\tn Q_\a\,H(1/\s_+)\,e_p.\]
Using the general fact $Q_\a\,W(v)\,(I-P_\a)\,e_p=e^{i\a p}\,P_\a\,H(\tl v)\,e_p$, which comes from
\[\int_\a^\iy V(\a-x-y)\,e^{ipy}\,dy=\int_0^\iy V(-x-y)\,e^{ip(\a+y)}\,dy,\]
we see that the above equals
\[{e^{i\a p}\ov c_+}\,Q_a\,H(\st)e_p\tn P_\a\,H(1/\smt)\,e_p+
{e^{i\a p}\ov c_-}\,P_\a\,H(\s)\,e_p\tn Q_\a\,H(1/\s_+)\,e_p\]
\[={e^{i\a p}\ov c_+}\,Q_a\,H(\st)e_p\tn P_\a\,\e+
{e^{i\a p}\ov c_-}\,P_\a\,H(\s)\,e_p\tn Q_\a\,\z,\]
in our notation.

This suggests that we add to our first guess a linear combination of tensor products with right factors $P_\a\,\e$ and $Q_\a\,\z$. By symmetry (or taking transposes of the above applied to $\st$) the left factors should be $P_\a\,\z$ and $Q_\a\,\e$. So it should be of the form
\[(a\,P_\a\,\z+b\,Q_\a\,\e)\tn\e+(c\,P_\a\,\z+d\,Q_\a\,\e)\tn Q_\a\,\z.\]
To see if this works, we compute $\W(\s)\,P_\a\,\z$ and $\W(\s)\,Q_\a\,\e$. We have
\[\W(\s)\,P_\a\,\z=P_\a\,W(\s)\,P_\a\,H(1/\s_+)\,e_p=P_\a\,W(\s)\,H(1/\s_+)\,e_p-P_\a\,W(\s)\,(I-P_\a)\,H(\s_+\inv)\,e_p\]
\be\eq-P_\a\,H(\s)\,W(1/\spt)\,e_p+{1\ov2ip\,c_+}P_\a\,W(\s)\,(I-P_\a)\,e_p\label{eq}\ee
\[=-{1\ov\s_+(p)}P_\a\,H(\s)\,e_p+{e^{i\a p}\ov2ip\,c_+}\,Q_\a\,H(\st)\,e_p.\]
Similarly
\[\W(\s)\,Q_\a\,\e=Q_\a\,\W(\st)\,\e\]
\[\eq-{1\ov\s_-(-p)}\,Q_\a\,H(\st)\,e_p+{e^{i\a p}\ov2ip\,c_-}\,P_\a\,H(\s)\,e_p.\]
This will cancel the previous extra terms when we add, if our coefficients are such that in $a\,P_\a\,\z+b\,Q_\a\,\e$ the coefficient of $P_\a\,H(\s)\,e_p$ is zero while the coefficient of $Q_\a\,H(\st)\,e_p$ is $-e^{i\a p}/c_+$, and in $c\,P_\a\,\z+d\,Q_\a\,\e$ the coefficient of $Q_\a\,H(\st)\,e_p$ is zero while that of $P_\a\,H(\s)\,e_p$ is $-e^{i\a p}/c_-$. We see that this does occur when
\[a={2ip\,\s_-(-p)\,\s_+(p)\,e^{2i\a p}\ov\D}={(2ip)^3\,\t_-(-p)\,\t_+(p)\,e^{2i\a p}\ov\D},\]
\[b={(2ip)^2\,c_-\,\s_-(-p)\,e^{i\a p}\ov\D}=-{(2ip)^3\,\t_-(p)\,\t_-(-p)\,e^{i\a p}\ov\D},\]
\[c={(2ip)^2\,c_+\,\s_+(p)\,e^{i\a p}\ov\D}=-{(2ip)^3\,\t_+(p)\,\t_+(-p)\,e^{i\a p}\ov\D},\]
\[d={2ip\,\s_-(-p)\,\s_+(p)\,e^{2i\a p}\ov\D}={(2ip)^3\,\t_-(-p)\,\t_+(p)\,e^{2i\a p}\ov\D},\]
where 
\[\D=(2ip)^2\,c_-\,c_+-e^{2i\a p}\,\s_-(-p)\,\s_+(p)={(2ip)^3\,
e^{i\a p}\ov A}.\]
Thus
\[a=\t_-(-p)\,\t_+(p)\,e^{i\a p}\,A,\ \ \ b=-\t_-(p)\,\t_-(-p)\,A,\]
\[c=-\t_+(p)\,\t_+(-p)\,A,\ \ \ d=\t_-(-p)\,\t_+(p)\,e^{i\a p}\,A.\]

This tells us that when we define $T_\a$ by (\ref{Ta}) then the first statement of the lemma holds. (The reason for the extra factor $A$ in the statement is that the error inherent in (\ref{eq}) gets multiplied by a factor proportional to $A$.) The second is obtained by applying the first to $\st$ and taking transposes. This completes the proof of the Lemma.
\sp 

\noi{\bf Proof of Theorem 1}. The theorem follows from the relations of the lemma by using the easily checked fact that $T_\a$ differs from $I$ by an integral operator whose kernel is $O(\a)$. 

\begin{center}{\bf 2. Proof of Theorem 2}\end{center}

Since the statement of the lemma to Theorem 1 holds uniformly in $p$, we need only compute the $p\to0$ limit of $T_\a$ to obtain its analogue when $p=0$. The computations simplify if we take $\t_\pm(0)=1$, which we may do since the theorem for the general case follows by applying the special case to the symbol $\s(\x)/\t(0)$. We easily compute that
\[A= -(\a+i(\t_-'(0)-\t_+'(0)))\inv+O(p)=-B\inv+O(p).\]

For the rest, a tedious but straightforward computation\footnote{The tedious computation was done by Maple.} shows that all terms of (\ref{Ta}) except for the first combine as $-B\inv$ times
\[[P_\a\,H(u_+)1-Q_\a\,H(\umt)1-x+i\t_+'(0)]\tn[P_\a\,H(\umt)1-Q_\a\,H(u_+)1-x-i\t_-'(0)]+O(p).\]
(All terms involving $p\inv$ cancel.) Now
\[W(1/\s_+)1=x+W(u_+)1,\]
while $H(u_+)1$ has Fourier transform
\[\({u_+-u_+(0)\ov0+i\x}\)_+=-\({u_+-u_+(0)\ov0-i\x}\)_+,\]
so
\[H(u_+)1=-W(u_+)1+u_+(0)=-W(u_+)1-i\t_+'(0).\]
(Observe that $\t_\pm(\x)=1\mp i\x\,u_\pm(\x)$.) Therefore
\[W(1/\s_+)1=x-H(u_+)1-i\t_+'(0),\]
and similarly for the right side of the tensor product. Therefore the tensor product equals
\[[P_\a\,W(1/\s_+)1+Q_\a\,H(\umt)1]\tn[P_\a\,W(1/\smt)1+Q_\a\,H(u_+)1].\]
Thus if we set
\[T_\a^0=P_\a\,W(1/\s_+)\,W(1/\s_-)\,P_\a
-Q_\a\,H(\wt{u_-})\,H(u_+)\,Q_\a\]
\be-B\inv\,[Q_\a\,H(\wt{u_-})1+P_\a\,W(1/\s_+)1]\tn[Q_\a\,H(u_+)1+P_\a\,W(1/\smt)1],\label{Ta0}\ee
then we have
\be\W(\s)\,T_\a^0=I+O(A\,e^{-\dl\a}),\quad T_\a^0\,\W(\s)=I+O(A\,e^{-\dl\a}).\label{Ta0W}\ee
This gives the statement of Theorem 2 when we check that $T_\a^0$ differs from $I$ by an integral operator whose kernel is $O(\a)$.
\sp

\begin{center}{\bf 3. Proof of Theorem 3}\end{center}

Suppose we want to prove (\ref{detpne0}) for a symbol $\s_1$. If we take any other $\s_0$ of the same form we can write $\s_1=\s_0\,e^{\r}$, and we define the family $\s=\s_\la=\s_0\,e^{\la\,\r}$ where $\la\in[0,\,1]$.  We assume first that $A$ satisfies the bound in the statement of Theorem~1 for all $\la\in[0,\,1]$, so that we can use its conclusion and the relation
\[{d\ov d\la}\log\,\d\W(\s)=\tr \W(\s)\inv\,\W(\s\r).\]

The operator on the right is
\[\W(\s)\inv\,[\W(\s)\,\W(\r)+P_\a\,H(\s)\,H(\rt)\,P_\a+
Q_\a\,H(\st)\,H(\r)\,Q_\a]\]
\be=\W(\r)+\W(\s)\inv\,[P_\a\,H(\s)\,H(\rt)\,P_\a+
Q_\a\,H(\st)\,H(\r)\,Q_\a].\label{PQ}\ee
We look first at 
\[\tr \W(\s)\inv\,P_\a\,H(\s)\,H(\rt)\,P_\a,\]
and in particular to the contribution to this of bracketed expression in (\ref{Ta}) with its factor $A$.

First, the factors $P_\a$ may be removed with error $\Oda$. Then $H(\rt)$ may be moved around to the left, so the trace of the product $\eq A$ times that of
\[e^{i\a p}\,\t_+(p)\,\t_-(-p)\,H(\rt)\,\z\tn H(\s)\,\e-\t_+(p)\,\t_+(-p)\,H(\rt)\,\z\tn H(\s)\,Q_\a\z\]
\[-\t_-(p)\,\t_-(-p)\,H(\rt)\,Q_\a \e\tn H(\s)\,\e+e^{i\a p}\,\t_+(p)\,\t_-(-p)\, H(\rt)\,Q_\a\z\tn H(\s)\, Q_\a\e.\]

The function $H(\s)\,Q_\a\z$ contains as a summand $H(\s)\,Q_\a\,H(u_+)\,e_p$ (the rest comes from the $e_p\tn e_p$ part oh $H(1/\s_+)$), and this is $\Oda$ since $H(\s)$ and $H(u_+)$ have exponentially decaying kernels. Similarly for $Q_\a\,\e$. Therefore with this error we may make the replacements
\[Q_\a\z\to -{e^{i\a p}\ov 2ip\,c_+}\,e_{-p}=-{e^{i\a p}\ov2ip\,\t_+(-p)}\,e_{-p},\]
\[Q_\a\e\to -{e^{i\a p}\ov 2ip c_-}\,e_{-p}=-{e^{i\a p}\ov2ip\,\t_+(p)}\,e_{-p}.\]

Next, for any plus function $v_+$ the Fourier transform of $H(v_+)e_p$ is
\[\({v_+\ov0+i(\x-p)}\)_+={v_+(p)-v_+\ov0-i(\x-p)}.\]
Therefore if we denote the Fourier transform by $\F$ we have in particular 
\[\F\,\e={\smt(p)\inv-\smt\inv\ov0-i(\x-p)}=-{(2ip\,\t_-(-p))\inv+\smt\inv\ov0-i(\x-p)},\]
and similarly
\[\F\,\z={\s_+(p)\inv-\s_+\inv\ov0-i(\x-p)}=-{(2ip\,\t_+(p))\inv+\s_+\inv\ov0-i(\x-p)}.\]

If we make these substitutions we find that we are left with four tensor products, only two of which have nonzero coefficients. We show them below with their coefficients.
 
\be H(\rt)\,e_{-p}\tn H(\s)\,e_{-p},\quad\quad -\,{e^{2i\a p}\ov 2ip\,\t_+(-p)\,\t_-(p)}\,A\inv,\label{prod1}\ee
\be H(\rt)\,\(\F\inv{\s_+\inv\ov0-i(\x-p)}\)\tn H(\s)\,\(\F\inv{\smt\inv\ov0-i(\x-p)}\),\quad\quad
e^{i\a p}\,\t_+(p)\,\t_-(-p).\label{prod2}\ee

Next we consider $P_\a\,H(\s)\,H(\rt)\,P_\a$ left-multiplied by the $Q_\a\,H(1/\smt)\,H(1/\s_+)\,Q_\a$ part of $\W(\s)\inv$. Again, because of the $Q_\a$ we may make the replacements
\[Q_\a\,H(1/\smt)\to c_-\inv\,e^{i\a p}e_{-p}\tn e_p,\quad
H(1/\s_+)\,Q_\a\to c_+\inv\,e^{i\a p}e_{p}\tn e_{-p},\]
the product of these being
\[-{e^{2i\a p}\ov 2ip\,c_-\,c_+}\,e_{-p}\tn e_{-p},\]
and after bringing the $H(\rt)$ around, the corresponding tensor product is $H(\rt)\,e_{-p}\tn H(\s)\,e_{-p}$ with the factor
\[-{e^{2i\a p}\ov 2ip\,\t_+(-p)\,\t_-(p)}.\]
This cancels the term in (\ref{prod1}) when we multiply it by $A$ and subtract this one.

To compute the trace of the term in (\ref{prod2}) we note that the factor on the right has Fourier transform
\[\(\s\,{\s_-\inv\ov 0+i(\x+p)}\)_+=
-{\s_+\ov0-i(\x+p)}\]
since $\s_+(-p)=0$.  Therefore the inner product of the two factors is 
\[-{1\ov2\pi}\int\r\,{\s_+\inv\ov0-i(\x-p)}\,{\s_+\ov0-i(\x+p)}\,d\x
=-{1\ov2\pi}\int{\r_-\ov(0-i(\x-p))\,(0-i(\x+p))}\,d\x\]
and so the trace with its factor (and the factor $A$) is 
\be A\,e^{i\a p}\,\t_+(p)\,\t_-(-p)\,{\r_-(p)-\r_-(-p)\ov2ip}.\label{P1}\ee

As for the $P_\a\,W(1/\s_+)\,W(1/\s_-)\,P_\a$ part of the inverse, this is $\eq$  the trace of
\[W(1/\s_+)\,W(1/\s_-)\,H(\s)\,H(\rt)=W(1/\s_+)\,H(\s_+)\,H(\rt)\]
\[=-H(1/\s_+)\,W(\spt)\,H(\rt)=-H(1/\s_+)\,H(\rt\,\spt).\]
The trace (including the minus sign) equals
\be{1\ov2\pi i}\int(\r\,\s_+)'\s_+\inv\,d\x={1\ov2\pi i}\int\r_-\,(\log \s_+)'\,d\x,\label{P2}\ee
since $\int \r'\,d\x=0$. 

So (\ref{P1})+(\ref{P2}) is the contribution of the $P_\a\,H(\s)\,H(\rt)\,P_\a$ term in (\ref{PQ}) to the trace.
For the contribution of the $Q_\a\,H(\st)\,H(\r)\,Q_\a$ term, observe that
\[\W(\s)\,\inv\,Q_\a\,H(\st)\,H(\r)\,Q_\a=
Q_\a\,\W(\st)\,\inv\,H(\st)\,H(\r)\,Q_\a,\]
which has the same trace as
\[\W(\st)\,\inv\,H(\st)\,H(\r)\,P_\a.\]
Hence we need only replace
$\s$ by $\st$ everywhere in what we derivated above. Then we add two and we get 
\pagebreak
\[A\,e^{i\a p}\,\t_+(p)\,\t_-(-p)\,{\r_-(p)-\r_-(-p)+\r_+(-p)-\r_+(p)\ov2ip}\]
\[+{1\ov2\pi i}\int[\r_-\,(\log \s_+)'-\r_+\,(\log \s_-)']\,d\x.\]

And finally we have to add to this
\[\tr\,\W(\r)={\a\ov2\pi}\int \r(\x)\,d\x.\]

What is this the derivative of? Since $(d/d\la)\log\s=\r$, for the last term the answer is
\[{\a\ov2\pi}\int \log\s(\x)\,d\x.\]
Here we use our chosen factorization of $\s$ to define the logarithm. For the next to last terms, since $(d/d\la)\log\s_\pm=\r_\pm$  we have 
\[{d\ov d\la}\int (\log\s_+)'\,\log\s_-\,d\x=\int [\r_+'\,\log\s_-+(\log\s_+)'\,\r_-]\,d\x,\]
so the middle term is $d/d\la$ of
\[{1\ov2\pi i}\int (\log\s_+)'\,\log\s_-\,d\x.\]

For the first term we compute that since $d\t_\pm/d\la=\r_\pm\,\t_\pm$,
\[{d\ov d\la}\log A={d\ov d\la}\log{1\ov \t_-(p)\,\t_+(-p)\,e^{-i\a p}-\t_+(p)\,\t_-(-p)\,e^{i\a p}}\]
\[={-(\r_-(p)+\r_+(-p))\,\t_-(p)\,\t_+(-p)\,e^{-i\a p}
+(\r_+(p)+\r_-(-p))\,\t_+(p)\,\t_-(-p)\,e^{i\a p}\ov 
\t_-(p)\,\t_+(-p)\,e^{-i\a p}-\t_+(p)\,\t_-(-p)\,e^{i\a p}}.\]
In the numerator we add $(\r_-(p)+\r_+(-p))\,\t_+(p)\,\t_-(-p)\,e^{i\a p}$ to the first summand and subtract it from the first. Then the quotient becomes
\[-\r_-(p)-\r_+(-p)-A\,e^{i\a p}\,\t_+(p)\,\t_-(-p)\,{\r_-(p)-\r_-(-p)+\r_+(-p)-\r_+(p)\ov2ip},\]
and so the first term is the logarithmic derivative of
\[(\t_+(-p)\,\t_-(p)\,A)\inv,\]
which is the same as the logarithmic derivative of
\[{\t_-(p)\,\t_+(-p)\,e^{-i\a p}-\t_+(p)\,\t_-(-p)\,e^{i\a p}\ov
\t_+(-p)\,\t_-(p)}.\]

Putting this together gives the logarithmic derivative of $\d\W(\s)$ with error $\Oda$. Integrating and exponentiating gives
\[{\d\W(\s)\ov G(\s)^\a}\eq C\,\exp\left\{{1\ov2\pi i}\int (\log\s_+)'\,\log\s_-\,d\x\right\}\]
\be\times {\t_-(p)\,\t_+(-p)\,e^{-i\a p}-\t_+(p)\,\t_-(-p)\,e^{i\a p}\ov\t_+(-p)\,\t_-(p)},\label{Cform}\ee
where $C$ is a constant independent of $\la$.

This was derived under the assumption that, for all $\la$,  $A=O(e^{\dl\a})$ for some small $\dl$, which means that $2i\a p$ is not exponentially close to any of the values of
\be\log{\t_-(p)\,\t_+(-p)\ov \t_+(p)\,\t_-(-p)},\label{tquotient}\ee
in particular if $2i\a p$ is bounded away from these values. Assume that this is so for our original $\s_1$, and suppose we had taken 
\[\s_0(\x)={\x^2-p^2\ov(q_1-i\x)\,(q_2+i\x)},\]
with $q_1,\,q_2$ in the right half-plane. The corresponding quotient in (\ref{tquotient}) is 
\[{q_1-ip\ov q_1+ip}\,{q_2-ip\ov q_2+ip}.\]
Each factor maps the right half-plane to the lower half-plane, so with proper choice of $q_i$ the ratio can be made arbitrarily close to the ratio in (\ref{tquotient}) coming from $\s_1$. The values of (\ref{tquotient}) for $\s_\la$ depend linearly on $\la$, so for this choice $2i\a p$ will be bounded away from the values of (\ref{tquotient}) for all $\la\in [0,\,1]$. Therefore (\ref{Cform}) holds for some $C$. We have \cite[Prop. 5.2]{B}
\[\d\W(\s_0)=e^{-\a(q_1+q_2)}\left[{(q_1+ip)\,(q_2+ip)\,e^{i\a p}-(q_1-ip)\,(q_2-ip)\,e^{-\a(p_1+p_2)/2}\ov2ip\,(q_1+q_2)}\right],\]
and a computation shows that for this $\s_0$ (\ref{Cform}) holds when $C=e^{i\a p}$.

This proves the theorem for our symbol $\s_1$ under the condition that $2i\a p$ is bounded away from the values of (\ref{tquotient}).
To remove the condition, suppose that it is very close to one of these values. Replace $\s_1$ by $\s_1\,e^{\la\r}$, with $\r$ now any symbol such that 
\[{\r_-(p)\,\r_+(-p)\ov \r_+(p)\,\r_-(-p)}\ne1,\]
and let $\la$ run over a little (but not too little) circle around zero. Then the corrrespoinding values of (\ref{tquotient}) will run over circles about the original ones, and so $2i\a p$ will be bounded away from them, uniformly in $\la$. By analyticity of the two sides in our formula, their values at $\la=0$  are obtained by integrating with respect to $\la$ over the circle. We deduce that the difference of the two sides, which is exponentially small uniformly for $\la$ on the circle, is also exponentially small for $\la=0$.

\begin{center}{\bf 4. Proof of Theorem 4}\end{center}

One can check that the right side of (\ref{detp=0}) is the limit of the right side of (\ref{detpne0}) as $p\to0$. Unfortunately it is not clear from our derivation that (\ref{detpne0}) holds uniformly in $p$, although it surely does. Fortunately the derivation of Theorem 4 from Theorem 2 is a lot easier than the derivation of Theorem 3 from Theorem 1, so we proceed with it.

We start in the same way, with (\ref{PQ}).
As before, the first term has trace
\[{\a\ov2\pi}\int\r(\x)\,d\x,\]
and we proceed to the contribution of 
\[\tr \W(\s)\inv\,P_\a\,H(\s)\,H(\rt)\,P_\a.\] 
The contribution of the $P_\a\,W(\s_+\inv)\,W(\s_-\inv)\,P_\a$ part of the inverse is 
\be{1\ov2\pi i}\int\r_-\,(\log \s_+)'\,d\x.\label{firstrs}\ee
as before. The contribution of the rest of $\W(\s)\inv$ is the trace of
\[-Q_\a\,H(\wt{u_-})\,H(u_+)\,Q_\a
-B\inv\,[Q_\a\,H(\wt{u_-})1+P_\a\,W(\s_+\inv)1]\tn[Q_\a\,H(u_+)1+P_\a\,W(\wt{\s_-}\inv)1]\]
right-multiplied by $P_\a\,H(\s)\,H(\rt)\,P_\a$. In computing the trace we may, as before, move the $H(\rt)\,P_\a$ around to the left. As we saw earlier, any trace involving two Hankel operator with nonsingular symbol with the operator $Q_\a$ between them will be $\Oda$. With similar error we may remove the factors $P_\a$. So with this error the trace reduces to that of
\[-B\inv\,H(\rt)\,W(\s_+\inv)1 \tn H(\s)\,W(\smt\inv)1.\]
The Fourier transform of the second factor is
\[\(\s\,{\s_-\inv\ov0+i\x}\)_+={\s_+\ov i\x},\]
while the Fourier transform of the first factor is
\[\(\rt\,{\spt\inv\ov0+i\x}\)_+.\]
The $\rt$ may be replaced by $\rt_+$, and the result after replacing $\x$ by $-\x$ is
\[\(\r_-\,{\s_+\inv\ov0-i\x}\)_-.\]
Notice that the expression in parenthesis has the singularity $(0-i\x)^{-2}$. 
When we multiply by the preceding and integrate we may remove the minus subscript, and we see that the trace in question equals 
\[{1\ov2\pi}\int{\r_-(\x)\ov(0-i\x)^2}\,d\x=i\,\r_-'(0).\]
So the contribution of the $P_\a\,H(\s)\,H(\rt)\,P_\a$ factor is
$-iB\inv\,\r_-'(0)$, 
which is to be added to (\ref{firstrs}). 
Similarly (or by applying the preceding to $\st$) we compute the contribution of the $Q_\a\,H(\st)\,H(\r)\,Q_\a$ term. Adding the two gives
\[{1\ov2\pi i}\int\r_-\,(\log \s_+)'\,d\x-{1\ov2\pi i}\int\r_+\,(\log \s_-)'\,d\x+iB\inv\,(\r_-'(0)-\r_+'(0)).\]
 
Thus the logarithmic derivative of the determinant equals
\[{\a\ov2\pi}\int\r(\x)\,d\x+{1\ov2\pi i}\int\r_-\,(\log \s_+)'\,d\x+{1\ov2\pi i}\int\r_+'\,\log \s_-\,d\x\]
\[+iB\inv\,(\r_-'(0)-\r_+'(0))+\Oda.\]
The first term is the logarithmic derivative of $G(\s)^\a$, and the next two terms combine as the logarithmic derivative of
\[\exp\left\{{1\ov2\pi i}\int (\log\s_+)'\,\log\s_-\,d\x\right\},\]
as we saw in the last section. Since (see footnote \ref{B})
\[B=\a-i\,\(\log{\s_+\ov\s_-}+\la\,(\r_+-\r_-)\)'(0),\]
the last term above is the logarithmic derivative of $B$.
This proves the theorem with a constant factor $C$ on the right side. Taking the special case $\s_0(\x)=\x^2/(1+\x^2)$ shows that $C=1$.
\pagebreak

\begin{center}{\bf 5. Comparison with the formula of \cite{AM}}\end{center}

We show here that (\ref{detpne0}) agrees with the result of \cite{AM}. Take $p>0$. We define two contours, both going from $-\iy$ to $\iy$. Contour $C_1$ has an indentation above $-p$ and below $p$, while $C_2$ has an indentation below $-p$ and above $p$. Each has its own geometric mean $G_i(\s)$ and constant $E_i(\s)$. These are defined by integrals over $C_i$ with integrands involving the associated factorizations. Our factorization $\s=\s_-\,\s_+$ is associated with $C_1$. Thus $G(\s)=G_1(\s),\ E(\s)=E_1(\s)$, and if we think of the right side of (\ref{detpne0}) as a sum then $G(\s)^\a$ times the first summand is exactly $G_1(\s)^\a\,E_1(\s)$.  

For the other terms we first express our exponential factor $E_1(\s)$ in terms of $E_2(\s)$. We have
\[\log E_1(\s)={1\ov2\pi i}\int_{C_1}\Big(\log((\x+p)\,\t_+)\Big)'\,\Big(\log((\x-p)\t_-)\Big)\,d\x\]
\[={1\ov2\pi i}\int_{-\iy-i\dl}^{\iy-i\dl}\Big(\log((\x+p)\,\t_+)\Big)'\,\Big(\log((\x-p)\t_-)\Big)\,d\x-\log(-2p\,\t_-(-p)),\]
\[\log E_2(\s)={1\ov2\pi i}\int_{C_2}\Big(\log((\x-p)\,\t_+)\Big)'\,\Big(\log((\x+p)\t_-)\Big)\,d\x\]
\[={1\ov2\pi i}\int_{-\iy-i\dl}^{\iy-i\dl}\Big(\log((\x-p)\,\t_+)\Big)'\,\Big(\log((\x+p)\t_-)\Big)\,d\x-\log(2p\,\t_-(p)),\]
since for the first we pass the pole at $\x=-p$ and for the second the pole at $\x=-p$. For the two resulting integrals, both over the contour in the lower half-plane, the first one minus the second equals
\[{1\ov2\pi i}\int(\log\t_+)'\,\log{\x-p\ov\x+p}\,d\x+
{1\ov2\pi i}\int\(\log{\x+p\ov\x-p}\)'\,\log((\x-p)\,\t_-)\,d\x\]
\[=
-{1\ov2\pi i}\int \log\t_+\,\(\log{\x-p\ov\x+p}\)'\,d\x+
{1\ov2\pi i}\int\(\log{\x+p\ov\x-p}\)'\,\log((\x-p)\,\t_-)\,d\x\]
\[=-\log\t_+(p)+\log\t_+(-p).\]
Putting these together shows that
\[E_1(\s)=-{\t_-(p)\ov\t_-(-p)}\,{\t_+(-p)\ov\t_+(p)}\,E_2(\s).\]

Finally, we compute that
$G_1(\s)/G_2(\s)=e^{-2ip}$.
Thus $e^{2i\a p}\,G(\s)^\a=G_2(\s)^\a$, so $G(\s)^\a$ times the second summand in our formula is exactly $G_2(\s)^\a\,E_2(\s)$. Hence our result can be stated
\[\d\W(\s)=G_1(\s)^\a\,(E_1(\s)+\Oda)+G_2(\s)^\a\,(E_2(\s)+\Oda),\]
which, except for the size of the error terms, agrees with \cite{AM}.

\begin{center}{\bf 6. Another proof of Theorem 2}\end{center}

As before, we may assume $\t_-(0)=\t_+(0)$. We begin by multiplying $\W(\s)$ on the right by $T_\a^0$, given by (\ref{Ta0}). For the product with $P_\a\,W(1/\s_+)\,W(1/\s_-)\,P_\a$ we obtain as before
\[\W(\s)\,P_\a\,W(1/\s_+)\,W(1/\s_-)\,P_\a=I-Q_a\,H(\st)\,H(1/\s_+)\,P_\a\,W(1/\smt)\,Q_\a,\]
but now we write
\[H(\st)\,H(1/\s_+)\,P_\a\,W(1/\smt)=H(\st)\,(1\tn1)\,P_\a\,W(1/\smt)+H(\st)\,H(u_+)\,P_\a\,W(\smt\inv)\]
\[\eq H(\st)\,(1\tn1)\,P_\a\,W(1/\smt)+H(\st)\,H(u_+)\,W(1/\smt)\]
\be=H(\st)\,(1\tn1)\,P_\a\,W(1/\smt)+H(\st)\,H(u_+/\s_-).\label{HHW}\ee
This is to be left- and right-multiplied by $Q_\a$.

Next, for the product with $Q_\a\,H(\umt)\,H(u_+)\,Q_\a$ we heve
\[\W(\s)\,Q_\a\,H(\umt)\,H(u_+)\,Q_\a=Q_\a\,W(\st)\,P_\a\,H(\umt)\,H(u_+)\,Q_\a\eq Q_\a\,W(\st)\,H(\umt)\,H(u_+)\,Q_\a,\]
and we compute
\[W(\st)\,H(\umt)\,H(u_+)=[H(\st\,\umt)-H(\st)\,W(u_-)]\,H(u_+)=H(\st\,\umt)\,H(u_+)-H(\st)\,H(u_-\,u_+)\]
\[=H\(\st\,\({1\ov\smt}-{1\ov0-i\x}\)\)\,H(u_+)-H(\st)\,H\(u_+\,\({1\ov\s_-}-{1\ov0+i\x}\)\)\]
\[=H\({\st\ov0+i\x}\)\,H(u_+)+H(\st)\,H\({u_+\ov0+i\x}\)-H(\st)\,H(u_+/\s_-).\]
Here we used the fact that $\st(0)=0$. The last term here cancels the same term appearing in (\ref{HHW}). The first two terms we may write as
\[H(\st)\,W\({1\ov0-i\x}\)\,H(u_+)+H(\st)\,W\({1\ov0+i\x}\)\,H(u_+).\]
Since
\[{1\ov0-i\x}+{1\ov0-i\x}=2\pi\dl=\hat{1},\]
these terms combine as
\[H(\st)\,(1\tn 1)\,H(u_+).\]

Putting together what we have done gives
\[\W(\s)\,[P_\a\,W(1/\s_+)\,W(1/\s_-)\,P_\a-Q_\a\,H(1/\smt)\,H(1/\s_+)\,Q_\a]\]
\[\eq I-Q_\a\,[H(\st)\,(1\tn1)\,P_\a\,W(1/\smt)+H(\st)\,(1\tn 1)\,H(u_+)]\,Q_\a.\]
For the right side of the first tensor product we use $Q_\a1=P_\a1$ to see that
\[(P_\a\,W(1/\smt)\,Q_\a)'1=Q_\a\,W(1/\s_-)\,P_\a\,1=P_\a\,W(1/\smt)\,Q_\a\,1\]
\[=P_\a\,W(1/\smt)\,P_\a\,1=P_\a\,W(1/\smt)\,1.\]
Thus the right side above equals
\[I-Q_\a\,H(\st)1\tn [P_\a\,W(1/\smt)1+Q_\a\,H(u_+)1].\]

Next, we will not make a guess as we did in the proof of the lemma since we really know the answer. So let us just continue by computing
\[\W(\s)\,[Q_\a\,H(\umt)1+P_\a\,W(\s_+\inv)1]=Q_\a\,W(\st)\,H(\umt)1+P_\a\,W(\s)\,P_\a\,W(\s_+\inv)\,P_\a\,1.\]
The first summand equals 
\[Q_\a\,H\(\st\,\({1\ov\smt}-{1\ov0-i\x}\)\)1-Q_\a\,H(\st)\,W(u_-)1=Q_\a\,H\({\st\ov0+i\x}\)\,1-u_-(0)\,Q_\a\,H(\st)\,1\]
\[=Q_\a\,H(\st)\,W\({1\ov0-i\x}\)1-u_-(0)\,Q_\a\,H(\st)\,1,\]

The second summand equals
\[P_\a\,W(\s_-)\,P_\a1-\a\,Q_\a\,H(\st)\,H(1/\s_+)\,Q_\a.\]
Using $H(1/\s_+)=1\tn 1+H(u_+)$ and $1\tn P_\a1=\a\,1\tn1$ we see that this is $\eq$ to
\[P_\a\,W(\s_-)\,P_\a1-\a\,Q_\a\,H(\st)1-Q_\a\,H(\st)\,H(u_+)1.\]
Since $H(u_+)1=-W(u_+)1+u_+(0)$ when we combine the terms we get
\[P_\a\,W(\s_-)\,P_\a1-\a\,Q_\a\,H(\st)1+Q_\a\,H(\st)\,W(\s_+\inv)1-(u_+(0)+u_-(0))\,Q_\a\,H(\st)1.\]
Since, as already observed, $u_\pm(0)=\mp i\,\t_\pm(0)$, the coefficient in the last term equals $-i\,(\t_+/\t_-)'(0)=i\,(\s_+/\s_-)'(0)$. Therefore the above equals
\[P_\a\,W(\s_-)\,P_\a1+Q_\a\,H(\st)\,W(\s_+\inv)\,1-B\,Q_\a\,H(\st)1.\]
Let us show that the first two terms cancel. We have $H(\st)\,W(\s_+\inv)\,1=H(\smt)\,1$. Since $\s_-(0)=0$ this has Fourier transform $\smt/i\x$, which is minus the Fourier transform of $W(\smt)\,1$. Therefore the second term above  equals
\[-Q_\a\,W(\smt)\,1=-Q_\a\,W(\smt)\,P_\a\,1=-P_\a\,W(\s_-)\,Q_\a\,1=-
P_\a\,W(\s_-)\,P_\a\,1.\]

We computed above what we obtain when we left-multiply the first two terms in (\ref{Ta0}) by $\W(\s)$. What we just computed shows that if we do the same with the last term and add, what results is $I+\Oda$. Similarly for multiplication on the right. (Or apply the preceding to $\st$ and take the transpose.) This gives (\ref{Ta0W}) and hence the statement of the theorem..

\begin{center}{\bf Acknowledgements}\end{center}

We thank Albrecht B\"ottcher for very helpful comments. This work was supported by National Science Foundation grant DMS-0552388.


\begin{thebibliography}{12}

\bibitem{AM} S. Albeverio and K. A. Makarov, {\it Extension of the Ahiezer-Kac determinnat formula to the case of real-valued symbols with two real zeros}, Acta Applic. Math. {\bf 62} (2000) 155--186.

\bibitem{B} A. B\"ottcher, {\it Wiener-Hopf determinants with rational symbols}, Math. Nachr. {\bf 144} (1989) 39--64.

\bibitem{BS1} A. B\"ottcher and B. Silbermann, {\it The asymptotic behavior of Toeplitz determinants for a generating function with zeros of integral order}, Math. Nachr. {\bf 102} (1981) 79--105.

\bibitem{BS} A. B\"ottcher and B. Silbermann, Analysis of Toeplitz Operators, Akademie-Verlag Berlin (1989).

\bibitem{BS2} A. B\"ottcher and B. Silbermann, Introduction to Large Truncated Toeplitz Matrices, Springer-Verlag New York (1998).

\bibitem{E} T. Ehrhardt, {\it A status report on the asymptotic behavior of Toeplitz determinants with Fisher-Hartwig singularities}, Oper. Th.: Adv. Appl. {\bf 124} (2001) 217--241.

\bibitem{M} L. V. Mikaelyan, {\it Asymptotics of determinants of truncated Wiener-Hopf operators in a singular case}, (Russian)  Akad. Nauk Armyan. SSR Dokl.  {\bf 82}  (1986) 151--155.

\bibitem{RS} P. Rambour and A. Seghier, {\it Exact and asymptotic inverse of the Toeplitz matrix with polynomial singular symbol}, C. R. Acad. Sci. Paris, Ser. I {\bf 335} (2002) 705--710.

\bibitem{W} H. Widom, {\it Extreme eigenvalues of N-dimensional convolution operators}, Trans. Amer.
Math. Soc. {\bf 106} (1963) 391--414.

\bibitem{W1} H. Widom {\it Asymptotics of a class of operator determinants with application to the cylindrical Toda equations}, arXiv: nlin.SI/0605021.

\bibitem{W2} H. Widom, {\it On the inverse and determinant of a truncated Wiener-Hopf operator}, arXiv: math.FA/0605076, v2.

\end{thebibliography}
\end{document}